\documentclass[a4paper,11pt,reqno]{amsart}

\usepackage[T1]{fontenc}
\usepackage{lmodern}
\usepackage{amsmath,amssymb,amsthm,mathtools}
\usepackage[top=3cm,bottom=3cm,left=3.15cm,right=3.15cm]{geometry}
\usepackage[colorlinks=true,linkcolor=blue,citecolor=blue,urlcolor=blue]{hyperref}
\usepackage{aliascnt}
\usepackage[capitalise,nameinlink,noabbrev]{cleveref}
\usepackage{microtype}
\usepackage{booktabs,array}
\usepackage{tikz}
\usepackage{tikz-cd}
\usepackage{needspace}

\AddToHook{env/proof/end}{\nopagebreak[4]}

\usepackage{etoolbox}
\makeatletter
\pretocmd{\@begintheorem}
  {\advance\@topsepadd 2.5pt plus 0.5pt minus 0.25pt\relax}
  {}{\PackageError{paragraph-spacing}{Theorem spacing patch failed}{}}
\expandafter\patchcmd\csname\string\proof\endcsname
  {\topsep6\p@\@plus6\p@\relax}
  {\topsep8.5\p@\@plus6.5\p@\@minus.25\p@\relax}
  {}{\PackageError{paragraph-spacing}{Proof spacing patch failed}{}}
\makeatother

\newtheorem{theorem}{Theorem}[section]
\newaliascnt{proposition}{theorem}
\newtheorem{proposition}[proposition]{Proposition}
\aliascntresetthe{proposition}
\newaliascnt{lemma}{theorem}
\newtheorem{lemma}[lemma]{Lemma}
\aliascntresetthe{lemma}
\newaliascnt{corollary}{theorem}
\newtheorem{corollary}[corollary]{Corollary}
\aliascntresetthe{corollary}
\newcommand{\restatedcorollaryname}{}
\newtheorem*{corollaryrestated}{\restatedcorollaryname}
\newenvironment{restatedcorollary}[1]
  {\renewcommand{\restatedcorollaryname}{Corollary~\ref*{#1}}%
   \begin{corollaryrestated}}
  {\end{corollaryrestated}}
\newcommand{\restatedtheoremname}{}
\newtheorem*{theoremrestated}{\restatedtheoremname}
\newenvironment{restatedtheorem}[1]
  {\renewcommand{\restatedtheoremname}{Theorem~\ref*{#1}}%
   \begin{theoremrestated}}
  {\end{theoremrestated}}
\theoremstyle{definition}
\newaliascnt{definition}{theorem}

\aliascntresetthe{definition}
\newtheoremstyle{spacedremark}
  {9pt plus 2pt minus 1pt}{9pt plus 2pt minus 1pt}
  {\normalfont}{}{\bfseries}{.}{0.5em}{}
\theoremstyle{spacedremark}
\newaliascnt{remark}{theorem}
\newtheorem{remark}[remark]{Remark}
\aliascntresetthe{remark}

\newcommand{\C}{\mathbb C}
\newcommand{\R}{\mathbb R}
\newcommand{\Z}{\mathbb Z}
\newcommand{\T}{\mathbb T}
\newcommand{\Pone}{\mathbb P^1}
\newcommand{\cO}{\mathcal O}
\newcommand{\Id}{\operatorname{Id}}

\newcommand{\HH}{\mathbb H}

\DeclareMathOperator{\diag}{diag}

\newcommand{\diff}{\mathrm{diff}}
\title[Complex structures on products with $S^4$]
{Complex Structures on Products with $S^4$}

\author[M. Gao]{Min Gao}
\address{Chern Institute of Mathematics and LPMC,
Nankai University, Tianjin 300071, P. R. China}
\email{mingao@mail.nankai.edu.cn}

\author[Z. Z. Tang]{Zizhou Tang}
\address{Chern Institute of Mathematics and LPMC,
Nankai University, Tianjin 300071, P. R. China}
\email{zztang@nankai.edu.cn}

\author[W. J. Yan]{Wenjiao Yan}
\address{School of Mathematical Sciences, Laboratory of Mathematics and Complex Systems,
Beijing Normal University, Beijing 100875, P. R. China}
\email{wjyan@bnu.edu.cn}

\date{}
\thanks{This work was partially supported by the NSFC
(Nos.~12271038, 12371048 and 12526205) and Nankai
Zhide Foundation.}

\hypersetup{
 pdftitle={Complex Structures on Products with S4},
 pdfsubject={Complex structures on four-sphere products from holomorphic line bundles},
 pdfkeywords={complex structures, equivariant splitting, torus bundles, nodal degeneration}}

\begin{document}
\begin{abstract}
The main theorem of the present paper is that $S^4\times P$ admits a complex structure whenever $P$ is
the fibre product, over a compact complex manifold, of the
unit circle bundles of an even number of holomorphic line
bundles.

This theorem has lots of significant applications. For example, $S^4\times G$ admits a complex structure for every compact Lie group $G$ of even dimension. The same
holds for $S^4\times V_{2s}(\C^n)$ and $S^4\times V_{2s}(\HH^n)$
when $n\geq 2s\geq2$, and for $S^4\times V_{4s}(\R^n)$ when
$n\ge4s\ge4$. The theorem also gives complex structures on the product
of $S^4$ with any even number of odd-dimensional spheres.
Moreover, complex structures exist on $S^4$ times some
moment-angle manifolds, including certain connected sums of sphere products.
The construction uses a holomorphic $(\C^*)^2$-action on $S^4 \times \T^2$ endowed with the complex structure we constructed in \cite{GaoTangYan2026} and a smooth splitting
equivariant under the compact torus.
\end{abstract}
\maketitle
\section{Introduction}
\label{sec:introduction}

Among spheres, only $S^2$ and $S^6$ admit almost complex structures \cite{BorelSerre1953}. In particular, $S^4$ admits no almost complex structure. 
Hopf's construction \cite{Hopf1948} equipped $S^1\times S^{2n+1}$ with a complex structure for every $n\geq1$. Calabi and Eckmann \cite{CalabiEckmann1953} later constructed complex structures on $S^{2p+1}\times S^{2q+1}$ for $p,q\geq1$. Bosio and Meersseman \cite{BosioMeersseman2006} constructed complex structures on certain connected sums of sphere products. Tang and Yan \cite{TangYan2022} constructed complex structures (integrable almost complex structures) on OT--FKM isoparametric hypersurfaces with $m=1$, including products of three spheres $S^1\times S^7\times S^6$ and $S^1\times S^3\times S^2$. Their complex structure on the latter manifold has vanishing first Chern class, while the product complex structure on a Hopf surface times $S^2$ has nonzero first Chern class. Qian, Tang and Yan \cite{QianTangYan2026} showed that the OT--FKM isoparametric hypersurfaces with $m=2$ also admit complex structures.

For products involving $S^4$, Calabi asked in 1958 whether
$S^4\times V$ admits a complex structure for a closed
orientable surface $V$ \cite{Calabi1958}. These products already admit almost complex structures, obtained
by embedding them as hypersurfaces in $\R^7$ and using the octonionic cross product.

A complex structure on a real smooth manifold $M$ of even dimension is insufficient in general for the
existence of even an almost complex structure on $S^4\times M$. For example, Tang
\cite{Tang1992} showed that $S^4\times\mathbb{CP}^2$ does not admit even an almost complex structure. Chakraborty and Thakur \cite{ChakrabortyThakur2016} extended the obstruction to $S^4\times\mathbb{CP}^{2n}$ for all $n\geq1$.
For products of $S^4$ with $S^{2n}$, only $S^4\times S^2$ admits an almost complex structure \cite{DattaSubramanian1990}. 
The existence of a complex structure in this remaining case is a separate question. Liang, Liu and Zhang \cite{LiangLiuZhang2026} have recently proposed a complex structure on $S^4 \times S^2$,  conditional on the validity of the torus-fibration construction of a complex structure on $S^6$.

In our preceding work \cite{GaoTangYan2026}, we constructed a complex threefold $X$ diffeomorphic to $S^4 \times \T^2$, answering
Calabi's question in genus one. In the present paper, we replace the product torus by the total space of a principal torus bundle determined by holomorphic line bundles.
For a holomorphic line bundle $L$ with a Hermitian metric, write $S(L)$ for its unit circle bundle and $L^\times$ for the complement of its zero section. We write $\T^r=\R^r/\Z^r$ and use the standard smooth structure on each sphere.

\begin{theorem}[Main theorem]\label{thm:main}
Let $L_1,\ldots,L_{2s}$ be holomorphic line bundles over a compact connected complex manifold $B$, where $s\geq1$. Denote the fibre product of their unit circle bundles over $B$ by
\[
 P=S(L_1)\times_B\cdots\times_B S(L_{2s}).
\]
Then $S^4\times P$ admits a complex structure.
\end{theorem}

We remark here that by \cite[Proposition~4.26]{FelixOpreaTanre2008}, the factor $P$ itself can be endowed with a complex structure, since the Chern connections on $L_j$ induce a connection on $P$ whose curvature is of type $(1,1)$.
The complex structures in \cref{thm:main} are necessarily
non-K\"ahler. Indeed, write $\dim_{\mathbb R}P=2m$. By the K\"unneth formula, every class in $H^2(S^4\times P;\mathbb R)$ is pulled back from $P$, so its $(m+2)$th power vanishes.
Thus $S^4\times P$ carries no symplectic form.

The odd-rank case follows from the main theorem by adding a trivial line bundle.

\begin{corollary}\label{cor:odd-rank}
If $P$ is a fibre product of the unit circle bundles of
an odd number of holomorphic line bundles over a compact
connected complex manifold, then $S^4\times P\times S^1$
admits a complex structure.
\end{corollary}

The following corollaries describe specific even-dimensional factors. The corresponding results for odd-dimensional factors follow from \cref{cor:odd-rank} by adding an $S^1$ factor. Accordingly, throughout the paper we state the odd-dimensional cases only when needed.

\begin{corollary}[Homogeneous spaces]
\label{cor:uniform-families}\label{cor:cohom-four-intro}
Let $G$ be a compact connected Lie group, $K$ the centralizer of a torus in $G$, and $L$ a closed normal subgroup of $K$. If $K/L\cong\T^{2s}$ for some $s\ge1$,
then $S^4\times G/L$ admits a complex structure.
\end{corollary}


We first recall Borel's results \cite{Borel1954}, which imply that the flag manifold $G/K$ carries a homogeneous K\"ahler structure. 
Over this flag base, $G/L$ can be realized as a fibre product of unit circle bundles of holomorphic line bundles. Then our main theorem applies.

In the case of compact connected Lie groups, Samelson \cite{Samelson1953} and Wang \cite{Wang1954} showed that every such group of even dimension carries a complex structure. We can prove that its product with $S^4$ also admits a complex structure whenever the group has even dimension.

\begin{corollary}[Compact Lie groups]\label{cor:intro-compact-groups}
For every compact connected Lie group $G$ of even
dimension, $S^4\times G$ admits a complex structure.
\end{corollary}

In particular, $S^4\times\mathrm{SU}(3)$ and $S^4\times G_2$ admit complex structures. However, neither $S^4 \times (\mathrm{SU}(3)/\T^2)$ nor $S^4 \times (G_2/\T^2)$ admits an almost complex structure; see \cref{hom:flag-product-obstructions}. Thus the positive rank assumption on $K/L$ in
\cref{cor:uniform-families} is essential.

For suitable choices of $G$, $K$, and $L$, the quotient $G/L$
is a Stiefel manifold, so \cref{cor:uniform-families} applies.
In \cref{hom:stiefel}, we also describe the corresponding
line bundles as tautological quotient lines over the flag base.

\begin{corollary}[Stiefel manifolds]\label{cor:intro-stiefel}
For $\mathbb F=\mathbb R$, $\mathbb C$, or $\mathbb H$, let
$V_j(\mathbb F^n)$ denote the Stiefel manifold of ordered
orthonormal $j$-frames in $\mathbb F^n$. For $s\ge1$, the manifolds
\[
 S^4\times V_{2s}(\C^n),\qquad
 S^4\times V_{2s}(\HH^n),\qquad
 S^4\times V_{4s}(\R^n)
\]
admit complex structures with $n\geq2s$ and $n>4s$, respectively. When $n=4s$, the last assertion
holds for the oriented component $\mathrm{SO}(4s)$.
\end{corollary}

Next, we write $L^{2a+1}(p)=S^{2a+1}/\Z_p$ for the lens space,
where $S^{2a+1}\subset\C^{a+1}$ is the unit sphere and
$\Z_p$ acts by scalar multiplication. Although products of scalar lens spaces are homogeneous, we obtain the following corollary directly from the unit circle bundles of the holomorphic line bundles $\mathcal O_{\mathbb P^a}(-p)$.

\begin{corollary}[Lens-space products]\label{cor:intro-lens-products}
Let $\ell\ge2$ be even, and let $a_j\ge0, p_j\ge1$
be integers for $1\leq j\leq\ell$. Then
\[
 S^4\times\prod_{j=1}^{\ell}L^{2a_j+1}(p_j)
\]
admits a complex structure.
\end{corollary}

In particular, taking all $p_j=1$, we obtain complex structures on
$S^4\times\prod_{j=1}^{2s}S^{2a_j+1}$ for all $s\ge1$
and $a_j\ge0$. 
These manifolds are stably parallelizable and have vanishing
Euler characteristic, so the existence of almost complex structures already follows from Sutherland's criterion
\cite{Sutherland1965}.
For $S^4\times S^5\times S^{2a+1}$, $a\geq0$,
complex structures can also be obtained by applying Morimoto's product theorem \cite[Theorems~2 and~3]{Morimoto1963} to the Sasakian structure on $S^4\times S^5$ constructed by Boyer and Galicki \cite[Theorem~2, $n=2$, $k=1$]{BoyerGalicki2006}.

It is quite meaningful that the corollaries above also apply to two families of isoparametric hypersurfaces. The quaternionic two-frame case of
\cref{cor:intro-stiefel} yields complex structures on the
definite OT--FKM isoparametric hypersurfaces $S^4\times V_2(\HH^k)$, $k\geq2$, answering the question of Qian, Tang and Yan \cite[Remark~2(ii)]{QianTangYan2026}.
The sphere-product case of \cref{cor:intro-lens-products}
includes the indefinite OT--FKM isoparametric hypersurfaces with multiplicities $(4,3)$, which are diffeomorphic to
$S^4\times S^3\times S^7$ \cite[Theorem~1.3]{QianTangYan2023}.

For simplicity, we let $\mathcal Z_\bullet$ denote the moment-angle manifold associated with the fan or polytope indicated by the subscript.

\begin{corollary}[Moment-angle manifolds]\label{cor:intro-moment-angle}
Let $\Sigma$ be a smooth complete rational fan in $\R^d$
with $m$ rays, where $m-d$ is positive and even.
Then $S^4\times\mathcal Z_\Sigma$ admits a complex structure,
with the smooth structure specified in \cref{tor:cox-construction}.
In particular, this applies to $S^4\times\mathcal Z_{P_m}$
for a Delzant $m$-gon $P_m$ with even $m\ge4$, and to
$S^4\times\mathcal Z_P$ when $P$ is obtained from a Delzant $d$-simplex by an odd number $q\ge1$ of vertex cuts, $d\ge2$;
see \eqref{tor:product_polygon} and \eqref{tor:product_truncated-simplex}.
\end{corollary}

The last case of \cref{cor:intro-moment-angle} states that the manifolds
\[
 S^4\times
 \left(
 \mathop{\#}_{j=1}^{q}
 j\binom{q+1}{j+1}
 \bigl(S^{j+2}\times S^{2d+q-j-1}\bigr)
 \right),
 \qquad d\ge2,\quad q\ge1\ \text{odd},
\]
admit complex structures.

For the smooth complete rational fans in
\cref{cor:intro-moment-angle}, the factors $\mathcal Z_\Sigma$
already carry complex structures by Panov--Ustinovsky
\cite[Theorem~3.3]{PanovUstinovsky2012}.
Our corollary gives complex structures on
$S^4\times\mathcal Z_\Sigma$.
Bosio--Meersseman \cite[Theorem~12.2]{BosioMeersseman2006}
construct complex structures on even-dimensional links
of the form $\T^\ell\times\mathcal Z_Q$ and on products with $S^1$ or diagonal circle quotients of odd-dimensional links, where $Q$ is a simple convex polytope and $\ell\ge0$.
We emphasize that the $S^4$ factor ensures that the homeomorphism types of the manifolds in \cref{cor:intro-moment-angle} do not occur in the Bosio--Meersseman construction; see \cref{tor:bm-examples}.

In \cref{sec:original-action,sec:equivariant}, we construct a
holomorphic action on $X$ and establish its equivariant splitting. \Cref{sec:main-proof} contains the proof of the main theorem. The applications are developed in
\cref{hom:section,hom:stiefel,ex:nonhomogeneous}.


\section{The original threefold and a holomorphic action}
\label{sec:original-action}

We recall the threefold $X$ constructed in \cite[Sections~3--6]{GaoTangYan2026} and equip it with a holomorphic action of $H\cong(\C^*)^2$. For two holomorphic line bundles $L_1 , L_2$ over $B$, $L_1^\times\times_B L_2^\times$ is a principal $H$-bundle. This action gives the holomorphic associated bundle used in \cref{sec:main-proof}. After reviewing the periods, nodal completions, and gluing of $X$, we extend the period translations across both singular fibres.

In additive notation,
\[
 H=(\C/\Z)^2, \qquad \T^2=(\R/\Z)^2\subset H.
\]
The exponential map
\[
 ([\alpha],[\beta])\longmapsto
 (e^{2\pi i\alpha},e^{2\pi i\beta})
\]
identifies $H$ with $(\C^*)^2$ and its subgroup $\T^2$
with $(S^1)^2$.

\subsection{Periods, completions, and gluing}

We begin with the regular torus family. For $t\in\C$, set \(q=e^{2\pi it}\in\C^*\) and \(B=\diag(0,1)\). The period matrix is
\begin{equation}\label{eq:periods}
 P(t)=P_{\text{tw}}(q)
 \begin{pmatrix}I_2&tB\\0&I_2\end{pmatrix},
\end{equation}
where
\[
P_{\text{tw}}(q)=\begin{pmatrix}1&-q&i&iq\\q&1&-iq&i\end{pmatrix}.
\]
By the construction in \cite[Section~3]{GaoTangYan2026}, the lattices \(\Lambda_t=P(t)\Z^4\) define a proper holomorphic family \(p:\mathcal G\to\C^*\) of complex two-tori. For a lift $t$ of $q$, we write a point as \((q,[\zeta])\), with \(\zeta\in\C^2\) taken modulo \(\Lambda_t\). The first and third columns are single-valued in $q$. They give the translation directions used below:
\begin{equation}\label{eq:invariant-columns}
 p_1(q)=(1,q)^{\mathsf T},\qquad
 p_3(q)=(i,-iq)^{\mathsf T}.
\end{equation}

For \(q\ne\pm i\), the normalized periods are given by
\begin{equation}\label{eq:normalized-periods}
 A(q)^{-1}P(t)=(I_2\;\;Z(t)),\qquad Z(t)=tB+C(q),
\end{equation}
where
\[
 A(q)=\begin{pmatrix}1&-q\\q&1\end{pmatrix},\qquad
 C(q)=\frac{i}{1+q^2}
 \begin{pmatrix}1-q^2&2q\\-2q&1-q^2\end{pmatrix}
 =\begin{pmatrix}a(q)&b(q)\\c(q)&d(q)\end{pmatrix}.
\]
The coefficients $a,b,c,d$ are holomorphic at $q=0$. Their values satisfy \(a(0)=d(0)=i\) and \(b(0)=c(0)=0\).
Choose $\epsilon>0$ sufficiently small and set $\Delta=\{|q|<\epsilon\}$. The normalized periods define a holomorphic family $\mathcal T^*$ over $\Delta^*=\Delta\setminus\{0\}$. The coordinate change $\widehat\zeta=A(q)^{-1}\zeta$ identifies the restriction of $\mathcal G$ over $\Delta^*$ biholomorphically with $\mathcal T^*$.

We next recall the completion at $q=0$ from \cite[Section~4]{GaoTangYan2026}. Its node coordinates will also be used to extend the action. After quotienting by the first two normalized period translations, introduce the multiplicative coordinates
\[
 z_1=e^{2\pi i\widehat\zeta_1},\qquad
 z_2=e^{2\pi i\widehat\zeta_2}.
\]
To add the central fibre, we use the chain \(\mathcal V=\coprod V_k/\sim\). Its node charts
\[
 V_k=\C_{z_1}^*\times
 \{(x_k,y_k,q)\in\C^2\times\Delta:x_ky_k=q\},\qquad
 k\in\Z,
\]
are glued by adjacent transitions
\begin{equation}\label{eq:node-transition}
 \chi_k:(z_1,x_k,y_k)_k\longmapsto
 (z_1,qx_k,x_k^{-1})_{k+1},\qquad x_k\ne0.
\end{equation}

Over $q\ne0$, the chain is identified with the multiplicative model above. Set \(u_f(q)=e^{2\pi if(q)}\) for $f=a,b,c,d$. The remaining period translations act on the chain by
\begin{equation}\label{eq:node-decks}
 \begin{aligned}
 \gamma_1(z_1,x_k,y_k)_k&=(u_a z_1,u_c^{-1}x_k,u_cy_k)_k,\\
 \gamma_2(z_1,x_k,y_k)_k&=(u_b z_1,u_d^{-1}x_k,u_dy_k)_{k+1}.
 \end{aligned}
\end{equation}
Let $\Gamma$ be the group generated by $\gamma_1$ and $\gamma_2$. As shown in \cite[Proposition~4.1]{GaoTangYan2026}, the quotient $\mathcal V/\Gamma$ is a complex threefold. It completes $\mathcal T^*$ across $q=0$, with
\[
(\mathcal V/\Gamma)|_{\Delta^*}\cong\mathcal T^*
\]
as holomorphic families over $\Delta^*$. The function $q=x_ky_k$ extends the projection to a proper holomorphic map $\pi:\mathcal V/\Gamma\to\Delta$. Its central fibre is
\[
 \pi^{-1}(0)\cong E_i\times C_{\mathrm{nod}},
\]
where $E_i=\C/(\Z+i\Z)$ and $C_{\mathrm{nod}}=\Pone/(0\sim\infty)$.

At $\infty$, set $s=-t$ and use the base coordinate $w=q^{-1}=e^{2\pi i s}$. On the punctured $w$-disc, the corresponding additive fibre coordinate is
\[
 \zeta_\infty=(w\zeta_2,-w\zeta_1)^{\mathsf T}.
\]
With respect to the local period matrix
\[
 P_\infty(s)=P_{\mathrm{tw}}(-w)
 \begin{pmatrix}I_2&sB\\0&I_2\end{pmatrix},
\]
the first two real period coordinates are unchanged and the last
two change sign:
\begin{equation}\label{eq:infinity-comparison}
 \zeta=P(t)(v_1,v_2,v_3,v_4)^{\mathsf T}
 \quad\Longrightarrow\quad
 \zeta_\infty=P_\infty(s)(v_1,v_2,-v_3,-v_4)^{\mathsf T}.
\end{equation}
In these coordinates, the normalized period matrix is
\[
 (I_2\;\;Z_\infty(s)),\qquad
 Z_\infty(s)=-Z(-s)=sB+C(-w).
\]
The same completion construction therefore applies at $\infty$. We replace $q$ by $w$ and $C(q)$ by $C(-w)$.
Gluing the local model $\mathcal V/\Gamma$ to $\mathcal G$ over the punctured $q$-disc gives $\overline{\mathcal G}_0$. The analogous construction at $\infty$ gives $\overline{\mathcal G}_\infty$. Write $\jmath_0$ and $\jmath_\infty$ for the inclusions of $\mathcal G$ into these completions.

Both completions contain the same regular family. To glue them, we use the translating section from \cite[Section~6]{GaoTangYan2026}. The corresponding fibrewise translation is $\tau_\sigma(q,[\zeta])=(q,[\zeta+u(t)])$, where
\begin{equation}\label{eq:affine-section}
 u(t)=P(t)\mathbf a(t),\qquad
 \mathbf a(t)=\bigl(0,-t(t+1)/2,0,t\bigr)^{\mathsf T}.
\end{equation}
The class of $u(t)$ modulo $\Lambda_t$ defines a single-valued holomorphic section over $\C^*$. We glue the two completions along their regular parts by the map
\begin{equation}\label{eq:holomorphic-gluing}
 \Phi:\jmath_0(\mathcal G)\longrightarrow\jmath_\infty(\mathcal G),
 \qquad
 \Phi=\jmath_\infty\circ\tau_\sigma\circ\jmath_0^{-1}.
\end{equation}
As proved in \cite[Section~6]{GaoTangYan2026}, the resulting space $X=\overline{\mathcal G}_0\bigcup_\Phi\overline{\mathcal G}_\infty$ is a compact complex threefold. It carries a proper holomorphic map $g:X\to\Pone$, with singular fibres only over $0$ and $\infty$.

\subsection{Extension of the period translations}

The periods $p_1,p_3$ define translations on $\mathcal G$. For the associated-bundle construction, we need these translations on all of $X$. We therefore extend them across both singular fibres.

\begin{proposition}\label{prop:holomorphic-action}
There is a holomorphic action
\(\rho:H\times X\to X\) preserving \(g\), given on the
regular family by
\begin{equation}\label{eq:regular-action}
 \rho_{\alpha,\beta}(q,[\zeta])
 =\bigl(q,[\zeta+\alpha p_1(q)+\beta p_3(q)]\bigr).
\end{equation}
\end{proposition}

\begin{proof}

On the regular family, changing $\alpha$ and $\beta$ by integers adds a lattice vector. Thus \eqref{eq:regular-action} is well defined for \(h=([\alpha],[\beta])\in H\). The columns $p_1$ and $p_3$ are single-valued in $q$, so the formula defines an action on all of $\mathcal G$.

We first extend the action at $q=0$. Fix representatives $\alpha,\beta\in\C$. In normalized coordinates, $p_1(q)$ and $p_3(q)$ become $(1,0)^{\mathsf T}$ and $(a(q),c(q))^{\mathsf T}$. The translation is therefore
\[
 \widehat\zeta\longmapsto
 \widehat\zeta+(\alpha+\beta a(q),\beta c(q))^{\mathsf T}.
\]
To check extension across the central fibre, use \(z_1=e^{2\pi i\widehat\zeta_1}\), \(y_k=q^{-k}e^{2\pi i\widehat\zeta_2}\), and \(x_k=q/y_k\). On the punctured part of a node chart, the translation becomes
\begin{equation}\label{eq:node-action}
 \rho_{\alpha,\beta}(z_1,x_k,y_k)_k
 =
 \left(
 e^{2\pi i(\alpha+\beta a(q))}z_1,\,
 e^{-2\pi i\beta c(q)}x_k,\,
 e^{2\pi i\beta c(q)}y_k
 \right)_k,\qquad q=x_ky_k.
\end{equation}
The coefficients $a,c$ are holomorphic at $q=0$, so the formula extends across the central fibre in every chart. Set \(u=e^{2\pi i(\alpha+\beta a(q))}\) and \(\mu=e^{2\pi i\beta c(q)}\). Compatibility with the chart transitions follows from
\[
 \begin{aligned}
 \chi_k\circ \rho_{\alpha,\beta}(z_1,x_k,y_k)_k
 &=(uz_1,q\mu^{-1}x_k,\mu x_k^{-1})_{k+1}\\
 &=\rho_{\alpha,\beta}\circ\chi_k(z_1,x_k,y_k)_k
 \qquad(x_k\ne0).
 \end{aligned}
 \]

Thus the local maps define transformations of the chain $\mathcal V$. We must also check that they descend through the period quotient. The generators $\gamma_1,\gamma_2$ satisfy
 \[
\begin{aligned}
\gamma_1\circ \rho_{\alpha,\beta}(z_1,x_k,y_k)_k
 &=(u_a u z_1,\,
     u_c^{-1}\mu^{-1}x_k,\,
     u_c\mu y_k)_k\\
 &=\rho_{\alpha,\beta}\circ\gamma_1(z_1,x_k,y_k)_k,\\[3pt]
 \gamma_2\circ \rho_{\alpha,\beta}(z_1,x_k,y_k)_k
 &=(u_b u z_1,\,
     u_d^{-1}\mu^{-1}x_k,\,
     u_d\mu y_k)_{k+1}\\
 &=\rho_{\alpha,\beta}\circ\gamma_2(z_1,x_k,y_k)_k.
\end{aligned}
\]

Hence the transformations descend to $\mathcal V/\Gamma$. It remains to check that they depend only on the classes of $\alpha$ and $\beta$ in $\C/\Z$.
On the chain $\mathcal{V}$, \(\rho_{1,0}=\Id\), \(\rho_{0,1}=\gamma_1\), so
\[
 \rho_{\alpha+m,\beta+n}=\gamma_1^n \rho_{\alpha,\beta},
 \qquad m,n\in\Z.
\]
The node formulas verify the group law and joint holomorphicity. The commutation relations and integer-shift identity above let these maps descend to a holomorphic $H$-action on $\mathcal V/\Gamma$. This action agrees with the original one over $q\ne0$.

At $\infty$, \eqref{eq:infinity-comparison} sends the translation direction \(\alpha p_1+\beta p_3\) to \(\alpha p_{1,\infty}-\beta p_{3,\infty}\). Here \(p_{j,\infty}\) is the \(j\)-th column of \(P_\infty\). The same calculation applies with $C(-w)$ in place of $C(q)$ and parameters $(\alpha,-\beta)$. This defines a holomorphic $H$-action on the completion at $\infty$.

It remains to check the gluing. Both $\rho_{\alpha,\beta}$ and $\tau_\sigma$ are fibrewise translations on $\mathcal G$, so they commute. Thus \eqref{eq:holomorphic-gluing} intertwines the actions on the two completed ends. These actions glue to a holomorphic $H$-action on $X$. Every local formula fixes the base coordinate, so the action preserves $g$.
\end{proof}

The action now supplies the holomorphic associated bundles in \cref{sec:main-proof}. We next prove a $\T^2$-equivariant smooth splitting of $X$ to determine their smooth type.

\section{Equivariant smooth splitting}\label{sec:equivariant}

We establish a $\T^2$-equivariant smooth splitting of $X$, which will determine the smooth structure of the associated bundle in \cref{sec:main-proof}. The first step is to make the relative trivializations of \cite[Section~5]{GaoTangYan2026} equivariant while retaining their prescribed boundary coordinates. A comparison of the two boundary maps then allows us to glue these trivializations. The remaining four-dimensional factor is the homotopy sphere identified in \cite[Section~7]{GaoTangYan2026}. Throughout this section, set
\[
 K_0=\{|q|\leq1\},\qquad K_\infty=\{|w|\leq1\},
 \qquad w=q^{-1},\qquad S=\{|q|=1\}.
\]

As in \cite[Section~2.2]{GaoTangYan2026}, let
\(N_0\to K_0\) and \(N_\infty\to K_\infty\) be the standard
one-node completions of the elliptic-curve families
\[
 \C/(\Z + (t+i)\Z),
 \qquad
 \C/(\Z + (s+i)\Z),
\]
over the punctured discs. The completions use local node models analogous to those in the preceding section, with the additional $z_1$ coordinate omitted. Each total space is a compact smooth four-manifold whose boundary lies over regular values.
On the regular part of \(N_0\), the coordinates are
\[
 q=re^{2\pi i\theta},\qquad
 t=\theta+\frac{\log r}{2\pi i},\qquad
 [\mu+(t+i)\nu]\in\C/(\Z+(t+i)\Z).
\]
These real coordinates \((r,[\theta,\mu,\nu])\) satisfy
\begin{equation}\label{eq:tate-seam}
 (\theta+1,\mu,\nu)\sim(\theta,\mu+\nu,\nu).
\end{equation}
For \(N_{\infty}\), the same coordinates and conventions hold with $w,s$ in place of $q,t$.

\subsection{An invariant relative trivialization}

The deformation in \cite[Section~5]{GaoTangYan2026} connects the end block to a product. We recall its parameter family and prescribed vector field before imposing equivariance. The periods are
\begin{equation}\label{eq:auxiliary-periods}
 P_\lambda(t)=P_{\mathrm{tw}}(\lambda q)\begin{pmatrix}I_2&tB\\0&I_2\end{pmatrix},
 \qquad0\leq\lambda\leq1,
\end{equation}
with \(P_1=P\) and
\[
 P_0(t)=\begin{pmatrix}1&0&i&0\\0&1&0&t+i\end{pmatrix}.
\]

Let \(\mathfrak B_I\) denote the smooth total space obtained by simultaneously completing this family, as in \cite[Section~5]{GaoTangYan2026}. It carries a proper map
\[
 \Pi=(\lambda,q):\mathfrak B_I\longrightarrow I\times\Delta_{1+\eta},
\]
where \(I\supset[0,1]\) is an open interval and \(\Delta_{1+\eta}\) is a disc, with \(\eta>0\). The slices
\(\mathcal B_\lambda=\Pi^{-1}(\{\lambda\}\times K_0)\) satisfy
\[
 \mathcal B_1=\overline{\mathcal G}_0|_{K_0},\qquad
 \mathcal B_0=E_i\times N_0.
\]
The construction supplies a smooth vector field \(V\) on
\(\mathfrak B_I\) satisfying
\begin{equation}\label{eq:parameter-field-conditions}
 d\lambda(V)=1,\qquad dq(V)=0.
\end{equation}
We use real period coordinates
$v=(v_1,v_2,v_3,v_4)^{\mathsf T}$, defined by
$[\zeta]=[P_\lambda(t)v]$. With $q=re^{2\pi i\theta}$,
these coordinates satisfy
\[
 (\theta+1,v_1,v_2,v_3,v_4)
 \sim(\theta,v_1,v_2+v_4,v_3,v_4).
\]
On a fixed open annulus about $S$, the vector field has the form
\[
 V=\left.\partial_\lambda\right|_{(r,\theta,v)},
\]
so its flow varies $\lambda$ while keeping $r$, $\theta$, and $v$ fixed.
This prescribed behavior is needed to compute the boundary gluing map exactly. The lemma below makes the flow equivariant by averaging $V$, without changing either \eqref{eq:parameter-field-conditions} or its annular expression. The underlying family and vector field are constructed in \cite[Section~5]{GaoTangYan2026}. For the associated-bundle construction in \cref{sec:main-proof}, we strengthen this relative trivialization to a $\T^2$-equivariant one.

\begin{lemma}\label{lem:equivariant-end}
There is a \(\T^2\)-equivariant diffeomorphism
\[
 D_0:\overline{\mathcal G}_0|_{K_0}\longrightarrow \T^2\times N_0.
\]
It extends across the cutting circle \(S\).
On an open annulus about that circle,
\begin{equation}\label{eq:equivariant-zero-marking}
 D_0(q,[P(t)v])
 =\bigl([v_1,v_3],(r,[\theta,v_2,v_4])\bigr).
\end{equation}
\end{lemma}

\begin{proof}

We first extend the torus action to the entire parameter family. On its regular part, define the action of
\(h=([\alpha],[\beta])\in\T^2\) by
\[
 \rho_{\alpha,\beta}(\lambda,q,[\zeta])
 =\bigl(\lambda,q,
 [\zeta+\alpha p_{1,\lambda}(q)+\beta p_{3,\lambda}(q)]\bigr),
\]
where \(p_{1,\lambda},p_{3,\lambda}\) are the first and third
columns of \(P_\lambda\), and \(\alpha,\beta\in\R\). The translations satisfy the group law and vary smoothly, so they define a smooth $\T^2$-action on the regular family.

The calculation in \cref{prop:holomorphic-action}, with
$C(q)$ replaced by $C(\lambda q)$, gives the node-chart formula:
\begin{equation}\label{eq:parameter-torus-action}
 \rho_{\alpha,\beta}(\lambda,z_1,x_k,y_k)_k
 =
 \left(\lambda,
 e^{2\pi i(\alpha+\beta a(\lambda q))}z_1,
 e^{-2\pi i\beta c(\lambda q)}x_k,
 e^{2\pi i\beta c(\lambda q)}y_k\right)_k.
\end{equation}
This formula extends smoothly across $q=0$ and preserves
$\lambda$ and $q=x_ky_k$. It also commutes with the chart
transitions and deck transformations of the parameter family:
\[
 \begin{aligned}
 \chi_k\circ\rho_{\alpha,\beta}
 &=\rho_{\alpha,\beta}\circ\chi_k,\\
 \rho_{\alpha,\beta}\circ\gamma_j
 &=\gamma_j\circ\rho_{\alpha,\beta},
 \qquad j=1,2.
 \end{aligned}
\]
The maps therefore descend to the completion. They depend only on $([\alpha],[\beta])$ and retain the group law, giving a smooth $\T^2$-action on $\mathfrak B_I$. Write $\rho_h=\rho_{\alpha,\beta}$.

The next step is to make the transport between parameter slices equivariant under this action. An invariant vector field has a local flow that commutes with the action. We therefore average $V$ using the normalized Haar measure $dh$ on $\T^2$:
\begin{equation}\label{eq:averaged-field}
 \overline V=\int_{\T^2}(\rho_h)_*V\,dh.
\end{equation}
The smoothness of the action and compactness of $\T^2$
ensure that $\overline V$ is smooth. Invariance of Haar measure
implies \((\rho_h)_*\overline V=\overline V\).
Since the action preserves \(\lambda\) and \(q\),
\[
 d\lambda(\overline V)=\int_{\T^2}1\,dh=1,\qquad
 dq(\overline V)=0.
\]

It remains to check that averaging preserves the annular marking. For $j=1,3$, we have $p_{j,\lambda}=P_\lambda(t)e_j$, so
\[
 [P_\lambda(t)v+\alpha p_{1,\lambda}+\beta p_{3,\lambda}]
 =[P_\lambda(t)(v+\alpha e_1+\beta e_3)].
\]
Thus, on the annulus,
\[
 \rho_{\alpha,\beta}(\lambda,r,\theta,v)
 =(\lambda,r,\theta,v+\alpha e_1+\beta e_3).
\]
This translation is independent of $\lambda$, hence
\[
 (\rho_h)_*
 \left(\left.\partial_\lambda\right|_{(r,\theta,v)}\right)
 =\left.\partial_\lambda\right|_{(r,\theta,v)}.
\]
Since the annulus is invariant and the prescribed field is unchanged by these translations, averaging leaves it fixed there:
\[
 \overline V=\int_{\T^2}V\,dh
 =V=\left.\partial_\lambda\right|_{(r,\theta,v)}.
\]

Let $\mathcal F_\tau$ be the flow of $\overline V$. Properness of $\Pi$ makes $\Pi^{-1}([0,1]\times K_0)$ compact. Since $\dot\lambda=1$ and $\dot q=0$, a trajectory starting in $\mathcal B_0$ remains in this set while it is defined for $0\leq\tau\leq1$. The continuation argument in \cite[Section~5]{GaoTangYan2026} gives the time-one map, and the reverse-time argument gives its inverse, so
\[
 \mathcal F_1:\mathcal B_0\longrightarrow\mathcal B_1
\]
is a diffeomorphism. Applying the same argument over a slightly larger closed disc extends it across the cutting circle.

The invariance $(\rho_h)_*\overline V=\overline V$ implies
\[
 \mathcal F_\tau\circ\rho_h=\rho_h\circ\mathcal F_\tau
\]
by \cite[Corollary~9.14]{Lee2013}. Thus $\mathcal F_1$ is
$\T^2$-equivariant. On the annulus,
\[
 \mathcal F_1(q,[P_0(t)v])=(q,[P_1(t)v]).
\]

At $\lambda=0$,
\[
 P_0(t)v=
 \begin{pmatrix}
 v_1+iv_3\\
 v_2+(t+i)v_4
 \end{pmatrix}.
\]
Thus $[v_1,v_3]$ are the real coordinates on $E_i$, and $[v_2,v_4]$ are the fibre coordinates on the regular fibres of $N_0$. The action translates $E_i$ by $[\alpha+i\beta]$ and fixes $N_0$. Composing $\mathcal F_1^{-1}$ with the identification
\[
 E_i\times N_0\longrightarrow\T^2\times N_0,\qquad
 ([x+iy],w)\longmapsto([x,y],w),
\]
we obtain the equivariant diffeomorphism $D_0$ with annular
formula \eqref{eq:equivariant-zero-marking}.
\end{proof}

We now apply the same construction at $\infty$, using the parameter family and prescribed vector field of \cite[Section~5]{GaoTangYan2026}, with periods
\[
 P_{\infty,\lambda}(s)
 =P_{\mathrm{tw}}(-\lambda w)\begin{pmatrix}I_2&sB\\0&I_2\end{pmatrix},\qquad s=-t.
\]
The proof of \cref{lem:equivariant-end} applies with \(\lambda q\) replaced by \(-\lambda w\). The change of coordinates \eqref{eq:infinity-comparison} introduces a reflection in the torus marking: the element
\(h=([\alpha],[\beta])\) acts by
\[
 \zeta_\infty\longmapsto
 \zeta_\infty+P_\infty(s)(\alpha,0,-\beta,0)^{\mathsf T}.
\]
The resulting diffeomorphism
\[
 D_\infty:\overline{\mathcal G}_\infty|_{K_\infty}
 \longrightarrow\T^2\times N_\infty
\]
satisfies
\[
 D_\infty(h\cdot p)=(R_T h)\cdot D_\infty(p),
 \qquad R_T[x,y]=[x,-y],
\]
where the target action is translation on the first factor. Thus $D_\infty$ is equivariant with the torus action modified by $R_T$. Its annular formula in the original coordinates is
\[
 D_\infty(q,[P(t)v])
 =\bigl([v_1,-v_3],(r^{-1},[-\theta,v_2,-v_4])\bigr).
\]
\subsection{The boundary identification and equivariant gluing}

The preserved annular markings allow us to compare the two end trivializations on the cutting circle. On $S$, we have $t=\theta\in\R$, and $\tau_\sigma$ adds the real period vector
\[
 \mathbf a(\theta)
 =\left(0,-\frac{\theta(\theta+1)}2,0,\theta\right)^{\mathsf T}.
\]
The annular formulas for $D_0$ and $D_\infty$ therefore give
\begin{equation}\label{eq:equivariant-boundary}
 f_S:=D_\infty\circ\Phi|_{g^{-1}(S)}\circ D_0^{-1}
 =R_T\times\psi:
 \T^2\times\partial N_0\longrightarrow\T^2\times\partial N_\infty,
\end{equation}
where
\[
 \psi[\theta,\mu,\nu]
 =\left[-\theta,\mu-\frac{\theta(\theta+1)}2,-\nu-\theta\right].
\]
The map $\psi$ is precisely the boundary diffeomorphism of
\cite[Section~7]{GaoTangYan2026}.

\begin{proposition}\label{prop:equivariant-product}
Let $Y=N_0\cup_\psi N_\infty$. There is a $\T^2$-equivariant diffeomorphism
\[
 X\cong_{\T^2}\T^2\times Y,
\]
where $\T^2$ acts on $X$ by \cref{prop:holomorphic-action}
and on $\T^2\times Y$ by translation
\[
 h\cdot(k,w)=(h+k,w), \qquad h,k\in \T^2, w \in Y. 
\]
\end{proposition}

\begin{proof}

Collar gluing gives $Y$ its smooth structure, as in \cite[Section~7]{GaoTangYan2026}. To give both end trivializations the same torus action, correct the marking at $\infty$ by setting
\[
 \widetilde D_\infty
 =(R_T\times\operatorname{id}_{N_\infty})\circ D_\infty.
\]
By \cref{lem:equivariant-end} and the corresponding formula
at $\infty$, for every \(h\in\T^2\), we have
\[
 \begin{aligned}
 D_0(h\cdot p)&=h\cdot D_0(p),\\
 \widetilde D_\infty(h\cdot p)
 &=(R_T\times\operatorname{id}_{N_{\infty}})
   \bigl((R_Th)\cdot D_\infty(p)\bigr)\\
 &=h\cdot\widetilde D_\infty(p).
 \end{aligned}
\]

Both block diffeomorphisms are now equivariant for the standard translation action. The same correction removes the torus reflection from the boundary map. Indeed, writing \(\Phi_S=\Phi|_{g^{-1}(S)}\), \eqref{eq:equivariant-boundary} gives
\[
 \begin{aligned}
 \widetilde D_\infty\circ\Phi_S
 &=(R_T\times\operatorname{id}_{N_{\infty}})\circ
   (R_T\times\psi)\circ D_0\\
 &=(\operatorname{id}_{\T^2}\times\psi)\circ D_0.
 \end{aligned}
\]
The block maps and their inverses therefore respect the boundary identifications. They induce an equivariant homeomorphism, smooth on each block,
\[
 D:X\longrightarrow\T^2\times Y.
\]

The actions are proper because $\T^2$ is compact. The equivariant homeomorphism $D$ restricts to a diffeomorphism on each block. By the equivariant smoothing theorem \cite[Theorem~9.4]{Kankaanrinta2007}, It can be adjusted in a $\T^2$-invariant neighbourhood of the common boundary to give
\[
 X\cong_{\T^2}\T^2\times Y.
\]
\end{proof}


The gluing map $\psi$ agrees with that in \cite[Section~7]{GaoTangYan2026}, so $Y$ is the same smooth four-manifold. In particular, $Y$ is a homotopy four-sphere by \cite[Proposition~7.2]{GaoTangYan2026}. This fact will be used with Hirsch's product theorem to determine the smooth type of the associated bundles in \cref{sec:main-proof}.

\section{Proof of the main theorem}
\label{sec:main-proof}

We now combine the holomorphic action and the equivariant smooth splitting of $X$ to prove the main theorem.
\begin{restatedtheorem}{thm:main}
Let $L_1,\ldots,L_{2s}$ be holomorphic line bundles over a compact connected complex manifold $B$, where $s\geq1$. Denote the fibre product of their unit circle bundles over $B$ by
\[
 P=S(L_1)\times_B\cdots\times_B S(L_{2s}).
\]
Then $S^4\times P$ admits a complex structure.
\end{restatedtheorem}

We first treat two line bundles. Motivated by the circle-bundle constructions of Sankaran and Thakur \cite{SankaranThakur2013}, we form a holomorphic associated bundle with fibre $X$. The compact torus reduction and the equivariant splitting identify its underlying smooth manifold with $P\times Y$. Hirsch's theorem then gives a diffeomorphism with $P\times S^4$. For more than two line bundles, the remaining circle factors provide a new complex base for the same construction.

\subsection{Two line bundles}

Set $H=(\C^*)^2$ and $\T^2=(S^1)^2$. We use the holomorphic $H$-action of \cref{prop:holomorphic-action} and the $\T^2$-equivariant smooth splitting of $X$ from \cref{prop:equivariant-product}. For associated bundles of a right principal $H$-bundle, our convention is
\[
 (qh,x)\sim(q,hx).
\]

\begin{proposition}\label{thm:rank-two}
Let $B$ be a compact connected complex manifold and let $L_1,L_2$
be holomorphic line bundles over $B$. Set
\[
 Q=L_1^\times\times_BL_2^\times,\qquad
 P=S(L_1)\times_BS(L_2).
\]
Then $\mathcal M=Q\times_HX$ is a compact complex manifold, and
\[
 \mathcal M\cong_{\diff}P\times S^4.
\]
\end{proposition}

\begin{proof}

First, the holomorphic principal bundle $Q\to B$ and the holomorphic $H$-action on $X$ give $\mathcal M$ a complex structure. In local holomorphic trivializations of $Q$, the associated bundle
\[
 \mathcal M=Q\times_HX\longrightarrow B
\]
has holomorphic transition maps and fibre $X$
\cite[Section~3.2(b)--(d)]{Hirzebruch1978}.
The compactness of both the base and the fibre makes $\mathcal M$ compact.

To identify the underlying smooth manifold, use the compact reduction $P\subset Q$ determined by the Hermitian metrics. The associated bundle $P\times_{\T^2}X\to B$ is smooth, and its total space is diffeomorphic to $\mathcal M$. Since every $q\in Q$ has a unique expression
\[
 q=pr,\qquad p\in P,\quad r\in(\R_{>0})^2,
\]
the map
\[
\begin{aligned}
     Q\times_HX &\longrightarrow P\times_{\T^2}X\\
    [q,x]&\longmapsto[p,rx]
\end{aligned}
\]
and its inverse
\[
[p,x]\longmapsto[p,x]
\]
are smooth inverse maps, giving
\[
 Q\times_HX\cong_{\diff}P\times_{\T^2}X.
\]
For well-definedness, replace $(q,x)$ by $(qh,h^{-1}x)$ and write $h=ta$, with $t\in\T^2$ and $a\in(\R_{>0})^2$. Its image is then $[pt,t^{-1}rx]=[p,rx]$.

The equivariant splitting now identifies the associated bundle with $P\times Y$. Let $\phi:\T^2\times Y\to X$ be the equivariant diffeomorphism of \cref{prop:equivariant-product}. The map
\[
\begin{aligned}
    P\times_{\T^2}X &\longrightarrow P\times Y\\
    [p,\phi(t,y)] &\longmapsto(pt,y)
\end{aligned}
\]
and the inverse map
\[
(p,y)\longmapsto[p,\phi(1,y)]
\]
are well defined because $\phi$ intertwines the torus actions. They are smooth inverses, so
\[
 P\times_{\T^2}X \cong_{\diff} P\times Y.
\]

Finally, $Y$ is a smooth homotopy four-sphere by \cite[Proposition~7.2]{GaoTangYan2026}. Since $P$ is closed and has positive dimension, Hirsch's theorem \cite[Theorem~6(a)]{Hirsch1965} gives
\[
P\times Y \cong_{\diff} P \times S^4.
\]

\Needspace{8\baselineskip}

These diffeomorphisms identify the compact complex manifold $\mathcal M$ with the required smooth product:
\[
\mathcal{M} \cong_{\diff} P \times S^4
\]
\end{proof}


\subsection{An even number of line bundles}

For $s>1$, we retain the first two line bundles and use the remaining $2s-2$ circle factors to form a new base. The Chern connections on the remaining line bundles have curvature of type $(1,1)$. Together with the even dimension of the torus fibre, this allows us to apply the result cited below. It gives the new base a complex structure with holomorphic projection to $B$. The holomorphic pullbacks of $L_1,L_2$ then reduce the proof to \cref{thm:rank-two}.
For $s>1$, let
\[
 \widehat B=S(L_3)\times_B\cdots\times_B S(L_{2s}),
\]
with natural projection $f:\widehat B\to B$.
By \cite[Proposition~4.26]{FelixOpreaTanre2008}, $\widehat B$ has the required complex structure and $f$ is holomorphic.
For $s=1$, set $\widehat B=B$ and $f=\operatorname{id}_B$.
Define
\[
 \widehat Q=(f^*L_1)^\times
 \times_{\widehat B}(f^*L_2)^\times,
 \qquad
 \widehat{\mathcal M}=\widehat Q\times_H X.
\]

\begin{proof}[Proof of \cref{thm:main}]

The case $s=1$ is \cref{thm:rank-two}. For $s>1$, the new base $\widehat B$ is compact and connected because its base and fibres are compact and connected. Since $f$ is holomorphic, so are $f^*L_1$ and $f^*L_2$. Thus the two-line-bundle construction applies over $\widehat B$. With the pullback metrics, the resulting unit-circle fibre product satisfies
\[
 \widehat P=S(f^*L_1)\times_{\widehat B}S(f^*L_2)
 \cong_{\diff}P.
\]
This identification restores the original ordering of the circle factors. It sends $(\hat b,v_1,v_2)$, with $\hat b=(b,v_3,\ldots,v_{2s})\in\widehat B$, to $(b,v_1,\ldots,v_{2s})$. Applying \cref{thm:rank-two} now gives the desired complex manifold:
\[
 \widehat {\mathcal M}= \widehat Q\times_HX
 \cong_{\diff}\widehat P\times S^4
 \cong_{\diff}P\times S^4.
\]
\end{proof}

For an odd number of line bundles, adjoining a trivial line bundle makes the number even and adds an $S^1$ factor to the total space.
\begin{restatedcorollary}{cor:odd-rank}
If $P$ is a fibre product of the unit circle bundles of
an odd number of holomorphic line bundles over a compact
connected complex manifold, then $S^4\times P\times S^1$
admits a complex structure.
\end{restatedcorollary}

\begin{proof}

Write the odd number as $2s+1$, with $s\geq0$. The trivial line bundle $\underline{\C}=B\times\C$ satisfies
\[
 P\times_B S(\underline{\C})\cong P\times S^1.
\]
Applying \cref{thm:main} to
$L_1,\ldots,L_{2s+1},\underline{\C}$ proves the assertion.
\end{proof}

With the existence proof complete, we record the holomorphic map inherited from $g:X\to\Pone$. The $H$-invariance of $g$ allows it to descend fibrewise to the associated bundle. The following proposition describes this map and its singular fibres.

\begin{proposition}\label{thm:degeneration}
The complex manifold $\widehat{\mathcal M}$ admits a proper
surjective holomorphic map
\[
 \widehat g:\widehat{\mathcal M}
 \longrightarrow\widehat B\times\Pone,
 \qquad
 [q,x]\longmapsto(\pi(q),g(x)),
\]
where $\pi:\widehat Q\to\widehat B$ is the bundle projection.
Its regular fibres are complex two-tori, and its singular
values are $\widehat B\times\{0,\infty\}$.
Each singular fibre is biholomorphic to
\[
 E_i\times C_{\mathrm{nod}},
 \qquad
 E_i=\C/(\Z+i\Z),\qquad
 C_{\mathrm{nod}}=\Pone/(0\sim\infty).
\]
\end{proposition}

\begin{proof}

The maps $g:X\to\Pone$ and $\pi:\widehat Q\to\widehat B$ are $H$-invariant and holomorphic, so
\[
 \widehat g([q,x]):=(\pi(q),g(x))
\]
is well defined and holomorphic. In a local holomorphic trivialization of $\widehat Q$, this map is the product of the identity on the base with $g$. Its surjectivity, fibres, and singular values therefore follow from those of $g$. Compactness of $\widehat{\mathcal M}$ gives properness.
\end{proof}


\section{Homogeneous spaces with flag bases}\label{hom:section}

We prove \cref{cor:uniform-families,cor:intro-compact-groups} by realizing the homogeneous spaces as unit-circle fibre products over flag manifolds. We first describe their smooth principal torus bundles, then use the flag bases to realize the associated line bundles holomorphically. This puts the spaces in the setting of \cref{thm:main}. All groups and quotients carry their standard smooth structures.

\subsection{Torus bundles over flag manifolds}\label{hom:flag-criterion}

\begin{restatedcorollary}{cor:uniform-families}
Let $G$ be a compact connected Lie group, $K$ the centralizer of a torus in $G$, and $L$ a closed normal subgroup of $K$. If $K/L\cong\T^{2s}$ for some $s\ge1$,
then $S^4\times G/L$ admits a complex structure.
\end{restatedcorollary}

\begin{proof}
We write $K=C_G(A)$ for a torus $A$ in $G$ and denote $F=G/K$.
There is a finite covering
\[
 Z(G)^0\times[G,G]\longrightarrow G,\qquad
 (z,h)\longmapsto zh,
\]
where $Z(G)^0$ is the identity component of the center
and $[G,G]$ is the commutator subgroup. $[G,G]$ is compact, connected, and semisimple.
Since $Z(G)^0 \subseteq K$ acts trivially on $F$,
we can reduce $G/K$ to the semisimple case.
Borel's results \cite[Theorems~2 and~4]{Borel1954} showed that $F$ is a simply connected homogeneous Kähler manifold, hence a Kähler $C$-space. Then by Griffiths' vanishing theorem \cite[Section~4, Theorem~$3'$]{Griffiths1963}, we obtain
\begin{equation}\label{hom:flag-vanishing}
 H^2(F,\mathcal O_F)=0.
\end{equation}
Here $\mathcal O_F$ is the sheaf of holomorphic functions. The vanishing \eqref{hom:flag-vanishing} will provide the holomorphic line bundles once we have described the smooth torus bundle by characters.

We next describe the smooth bundle over $F$. Denote $P=G/L$. Since $L\triangleleft K$, the projection
\[
 \pi:P\longrightarrow F,\qquad gL\longmapsto gK,
\]
has a well-defined right $K/L$-action
\[
 (gL)\cdot(kL)=gkL.
\]
A local smooth section $\sigma:U\to G$ of $G\to F$ defines
an equivariant trivialization
\[
 U\times K/L\longrightarrow\pi^{-1}(U),
 \qquad (b,kL)\longmapsto\sigma(b)kL.
\]
Thus $P\to F$ is a smooth principal $K/L$-bundle.

To express this principal bundle as a unit-circle fibre product,
we fix an identification $K/L\cong(S^1)^{2s}$ and take the
coordinate characters
\[
 \chi_j:K/L\longrightarrow S^1,\qquad
 \chi_j(t_1,\ldots,t_{2s})=t_j,
 \qquad 1\le j\le 2s.
\]
These characters form a basis of the character lattice of $K/L$.
The associated smooth complex line bundles are defined by
\[
 E_j=P\times_{\chi_j}\C,\qquad
 (pt,z)\sim(p,\chi_j(t)z).
\]
Since each $\chi_j$ takes values in $S^1$, 
$|[p,z]|_{E_j}=|z|$ is a well-defined Hermitian metric on $E_j$.
Consider the map
\[
 P\longrightarrow S(E_1)\times_F\cdots\times_F S(E_{2s}),
 \qquad p\longmapsto([p,1],\ldots,[p,1]).
\]
In the trivializations induced by a local section of $P$,
this map is
\[
 (b,t)\longmapsto
 \bigl(b,\chi_1(t),\ldots,\chi_{2s}(t)\bigr).
\]
Since $(\chi_1,\ldots,\chi_{2s})$ is the chosen identification of $K/L$ with $(S^1)^{2s}$, this is an isomorphism of smooth principal torus bundles.

It remains to realize the $E_j$ as holomorphic line bundles. This is where \eqref{hom:flag-vanishing} enters. The exponential sequence of sheaves
\[
 0\longrightarrow\underline{\Z}
 \longrightarrow\cO_F
 \xrightarrow{\,f\mapsto e^{2\pi i f}\,}
 \cO_F^*
 \longrightarrow 1,
\]
induces the exact segment
\[
 \operatorname{Pic}(F)=H^1(F,\cO_F^*)
 \xrightarrow{c_1}H^2(F;\Z)
 \longrightarrow H^2(F,\cO_F),
\]
where $\cO_F^*$ is the sheaf of nowhere-vanishing holomorphic
functions and $\underline{\Z}$ is the constant sheaf of integers.
Here we identify the cohomology of the constant sheaf
$\underline{\Z}$ with singular cohomology, as $F$ is a smooth
manifold.
Since \(H^2(F,\mathcal{O}_F) = 0\), the coboundary map \(c_1\) is surjective. \(c_1(E_j)\) is therefore the first Chern class of a holomorphic line bundle \(L_j\).
Smooth complex line bundles are classified by their first Chern
class, so there are smooth line-bundle isomorphisms
\(\phi_j:E_j\cong L_j\). We equip $L_j$ with the metric transported by $\phi_j$. Therefore
\begin{equation}\label{hom:unit-frame-map}
 \begin{aligned}
 P&\longrightarrow
 S(L_1)\times_F\cdots\times_F S(L_{2s}),\\
 p&\longmapsto
 \bigl(\phi_1([p,1]),\ldots,\phi_{2s}([p,1])\bigr)
 \end{aligned}
\end{equation}
is a diffeomorphism of principal torus bundles.

Thus $G/L$ is a unit-circle fibre product of $2s$ holomorphic line bundles over the compact connected complex manifold $F$. Since $s\geq1$, \cref{thm:main} applies and gives the conclusion.
\end{proof}

\subsection{Compact Lie groups}

The maximal torus bundle $\T_G\to G\to G/\T_G$ has the compact Lie group itself as total space. We apply \cref{cor:uniform-families} to this bundle. For $G$ of positive even dimension, the remaining point is to verify that its maximal torus also has even dimension.

\begin{restatedcorollary}{cor:intro-compact-groups}
For every compact connected Lie group $G$ of even
dimension, $S^4\times G$ admits a complex structure.
\end{restatedcorollary}

\begin{proof}
For a compact connected Lie group \(G\), we choose a  maximal torus \(\T_G\). Its centralizer in \(G\) is
\(K = \T_G\) \cite[Chapter~IV, Theorem~2.3(i)]{BrockerTomDieck1985}, and the flag base is \(G/\T_G\). Let \(L = \{e\}\), then
\[
K/L=\T_G,\qquad G/L=G.
\]
Thus the relevant bundle in \cref{cor:uniform-families} is simply \(\T_G\to G\to G/\T_G\). We now check the positivity and parity of its torus rank.

Write $r=\dim\T_G$.
Since $C_G(\T_G)=\T_G$, a trivial maximal torus would
force $G=\{e\}$. Thus $r>0$.
The flag manifold $G/\T_G$ has a standard complex structure
\cite[Theorem~2]{Borel1954}, so
\[
 \dim G-r=2\dim_{\C}(G/\T_G).
\]
Hence $\dim G\equiv r\pmod2$.
Under the hypothesis of the corollary, $r$ is positive and even.
\cref{cor:uniform-families} applies with
$K=\T_G$ and $L=\{e\}$.
\end{proof}

\begin{remark}\label{hom:flag-product-obstructions}
For $K=\T_G$, the choices $L=\{e\}$ and $L=K$ correspond
to $G$ and $G/\T_G$, respectively, with
$\dim(K/L)=0$ in the latter case.
However, neither $S^4\times(\mathrm{SU}(3)/\T^2)$ nor
$S^4\times(G_2/\T^2)$ admits an almost complex structure,
so the positive rank assumption in
\cref{cor:uniform-families} is essential.

The following sketch explains the obstruction in these two examples with $\dim(K/L)=0$.
Write \(M_1=S^4\times(\operatorname{SU}(3)/\T^2)\) and
\(M_2=S^4\times(G_2/\T^2)\).
If \(M_i\) admitted an almost complex structure \(J_i\),
the complex vector bundle \(E_i=(TM_i,J_i)\) would satisfy
\(\operatorname{ch}_{2j}(E_i)=0\) for \(j\ge1\),
since \(M_i\) is stably parallelizable \cite[Section~3.1]{TangXieYan2012}.
Then, using Bott periodicity, line-bundle generators, and the
cohomology rings; see
\cite{Atiyah1967,AtiyahHirzebruch1961,BurfittGrbic2023},
the calculation shows that
\[
 c_5(E_1)\in8H^{10}(M_1;\Z),\qquad
 c_8(E_2)\in16H^{16}(M_2;\Z).
\]
This contradicts
\(\langle c_{\mathrm{top}}(E_i),[M_i]\rangle=\chi(M_i)\),
since \(\chi(M_1)=12\) and \(\chi(M_2)=24\)
\cite[Section~14 and Corollary~11.12]{MilnorStasheff1974}.
\end{remark}

\section{Stiefel manifolds}\label{hom:stiefel}

We apply the homogeneous-space criterion to prove \cref{cor:intro-stiefel} for the complex, quaternionic, and real Stiefel manifolds. The proof identifies each frame space as $G/L$ and checks the rank of the torus $K/L$. We then describe the corresponding holomorphic line bundles over the flag bases using tautological quotients.

\begin{restatedcorollary}{cor:intro-stiefel}
For $\mathbb F=\mathbb R$, $\mathbb C$, or $\mathbb H$, let
$V_j(\mathbb F^n)$ denote the Stiefel manifold of ordered
orthonormal $j$-frames in $\mathbb F^n$. For $s\ge1$, the manifolds
\[
 S^4\times V_{2s}(\C^n),\qquad
 S^4\times V_{2s}(\HH^n),\qquad
 S^4\times V_{4s}(\R^n)
\]
admit complex structures with $n\geq2s$ and $n>4s$, respectively. When $n=4s$, the last assertion
holds for the oriented component $\mathrm{SO}(4s)$.
\end{restatedcorollary}

\subsection{Homogeneous presentations and proof}

\begin{proof}[Proof of \cref{cor:intro-stiefel}]

We treat complex and quaternionic $r$-frames, with $1\le r\le n$, and real $2t$-frames, with $1\le2t\le n$. To apply the criterion, we choose the groups \(G\), tori \(A\), centralizers $K=C_G(A)$, and subgroups $L$ in the following table:
\begin{center}
\begingroup
\setlength{\tabcolsep}{5pt}
\renewcommand{\arraystretch}{1.2}
\begin{tabular*}{\linewidth}{@{\extracolsep{\fill}}cccc@{}}
\toprule
$G$ & $A$ & $K=C_G(A)$ & $L$ \\
\midrule
$\operatorname{U}(n)$
& $\left\{
   \begin{gathered}
   \diag(z_1,\ldots,z_r,I_{n-r}):
   z_j\in S^1
   \end{gathered}
   \right\}$
& $\begin{gathered}
   \operatorname{U}(1)^r
   {}\times\operatorname{U}(n-r)
   \end{gathered}$
& $\operatorname{U}(n-r)$ \\
\addlinespace[8pt]
$\operatorname{Sp}(n)$
& $\left\{
   \begin{gathered}
   \diag(z_1,\ldots,z_r,I_{n-r}):
   z_j\in S^1
   \end{gathered}
   \right\}$
& $\begin{gathered}
   \operatorname{U}(1)^r
   {}\times\operatorname{Sp}(n-r)
   \end{gathered}$
& $\operatorname{Sp}(n-r)$ \\
\addlinespace[8pt]
$\operatorname{SO}(n)$
& $\left\{
   \begin{gathered}
   \diag(R(\theta_1),\ldots,R(\theta_t),I_{n-2t}):
   \theta_j\in\R
   \end{gathered}
   \right\}$
& $\begin{gathered}
   \operatorname{SO}(2)^t
   {}\times\operatorname{SO}(n-2t)
   \end{gathered}$
& $\operatorname{SO}(n-2t)$ \\
\bottomrule
\end{tabular*}
\endgroup
\end{center}
where 
\[
 R(\theta)=\begin{pmatrix}\cos\theta&-\sin\theta\\
                         \sin\theta&\cos\theta\end{pmatrix}.
\]
In the quaternionic case, $S^1$ lies in a fixed copy of
$\C\subset\HH$. The standard matrix groups $\operatorname{U}(n)$ and
$\operatorname{Sp}(n)$ preserve the standard complex and
quaternionic Hermitian forms, respectively, while
$\operatorname{SO}(n)$ preserves the Euclidean inner product
and orientation.

The torus $A$ rotates the first $r$ coordinate lines, or the
first $t$ coordinate two-planes, independently and fixes
the remaining subspace. Any matrix commuting with $A$
must preserve each of these subspaces.
On each complex or quaternionic coordinate line, it must act by multiplication by a unit complex number. On each real coordinate two-plane, it must act by a rotation.
These observations determine the centralizers in the table.
In each row, $L$ is the kernel of the projection of $K$ onto its circle factors, hence a closed normal subgroup. Thus $K/L\cong\T^r$ in the first two rows and $\T^t$ in the third. These are the torus fibres of $G/L\to G/K$.

Every orthonormal frame can be extended to an orthonormal basis. In the real case, a positively oriented extension exists when $n>2t$. For $n=2t$, we consider positively oriented full frames. Thus $G$ acts transitively on the corresponding frame space, with the stabilizer $L$ at the standard frame.
The orbit map therefore identifies $G/L$ with the following frame spaces:
\[
 G/L\cong_{\diff}
 \begin{cases}
 V_r(\C^n),&G=\operatorname{U}(n),\\
 V_r(\HH^n),&G=\operatorname{Sp}(n),\\
 V_{2t}(\R^n),&G=\operatorname{SO}(n),\ n>2t.
 \end{cases}
\]
For $n=2t$, $G/L \cong_{\diff} \operatorname{SO}(2t)$ instead of $V_{2t}(\R^{2t}) = \operatorname{O}(2t)$.

For complex and quaternionic frames, $K/L\cong\T^r$; for real $2t$-frames, $K/L\cong\T^t$. Taking $r=2s$ in the first two cases and $t=2s$ in the third makes the torus rank positive and even. The identifications above therefore allow us to apply \cref{cor:uniform-families}, proving all assertions.
\end{proof}

For $n=2t=4s$, the full frame space $V_{4s}(\R^{4s})$ is diffeomorphic to $\mathrm O(4s)$. Although \cref{cor:intro-stiefel} does not apply to it,  \cref{cor:intro-compact-groups}, applied componentwise,  shows that $S^4\times V_{4s}(\R^{4s})$
also admits a complex structure.

The quaternionic two-frame case also applies to the definite OT--FKM isoparametric hypersurfaces with $m=4$ and $l=4k$, $k\geq2$, since
\[
 M\cong_{\diff}S^4\times M_+\cong_{\diff}S^4\times V_2(\HH^k).
\]
The first diffeomorphism is recalled in \cite[Remark~2(ii)]{QianTangYan2026}, and the second follows from \cite{QianTangYan2022}.
Thus \cref{cor:intro-stiefel} implies that these
hypersurfaces admit complex structures.

\subsection{Tautological quotient line bundles}

We now realize the character line bundles in the proof of \cref{cor:uniform-families} explicitly. The flag bases $G/K$ have concrete complex models on which successive tautological quotients are holomorphic line bundles. Their unit-circle fibre products will recover $G/L$. We begin with the following three flag spaces:
\[
 \begin{aligned}
 F^{\C}_{n,r}
 &=\{V_\bullet=(0=V_0\subset V_1\subset\cdots\subset V_r
       \subseteq\C^n):\dim_{\C}V_j=j\},\\
 F^{\HH}_{n,r}
 &=\{V_\bullet=(0=V_0\subset V_1\subset\cdots\subset V_r
       \subseteq\C^{2n}):\dim_{\C}V_j=j,\ \omega|_{V_r}=0\},\\
 F^{\R}_{n,t}
 &=\{V_\bullet=(0=V_0\subset V_1\subset\cdots\subset V_t
       \subseteq\C^n):\dim_{\C}V_j=j,\ b|_{V_t}=0\}^{\circ}.
 \end{aligned}
\]
Here $\omega$ and $b$ are the standard complex symplectic form and symmetric bilinear form, respectively, given by
\[
 \omega\bigl((z,w),(z',w')\bigr)
 =\sum_{a=1}^n(z_aw'_a-w_az'_a),
 \qquad
 b(u,v)=\sum_{a=1}^n u_av_a.
\]
The standard action of $\operatorname{Sp}(n)$ on $\C^{2n}$ preserves $\omega$ and the standard action of $\operatorname{SO}(n)$ on $\C^n$ preserves $b$.

We write $F$ for the corresponding flag space and use
\[
 (E,\ell)=(\C^n,r),\quad(\C^{2n},r),\quad(\C^n,t),
\]
respectively. The superscript $\circ$ selects the connected
component containing the standard flag
\[
 V_j^0=\operatorname{span}_{\C}\{a_1,\ldots,a_j\},
 \qquad
 a_j=
 \begin{cases}
 e_j,&F=F^{\C}_{n,r}\text{ or }F^{\HH}_{n,r},\\
 (e_{2j-1}-i e_{2j})/\sqrt2,&F=F^{\R}_{n,t}.
 \end{cases}
\]
Here $e_j$ denotes the $j$th standard basis vector in $\C^n$ or $\C^{2n}$, as appropriate.

These spaces are smooth closed subvarieties of $\prod_{j=1}^{\ell}\operatorname{Gr}_j(E)$ and carry the induced complex structures. They model the flag bases: the action of $G$ on $F$ is transitive, with stabilizer $K$ at the standard flag $V_\bullet^0$. Hence the orbit map induces a diffeomorphism
\[
 G/K\longrightarrow F,\qquad gK\longmapsto gV_\bullet^0.
\]

The coordinate projections
$F\to\operatorname{Gr}_j(E)$, $V_\bullet\mapsto V_j$,
are holomorphic. The tautological bundle over $\operatorname{Gr}_j(E)$ is
\[
 \mathcal U_j
 =\{(W,v)\in\operatorname{Gr}_j(E)\times E:\ v\in W\},
 \qquad (W,v)\longmapsto W.
\]
Its pullback under $V_\bullet\mapsto V_j$ is therefore
\[
 \mathcal S_j
 =\{(V_\bullet,v)\in F\times E:\ v\in V_j\},
 \qquad (V_\bullet,v) \longmapsto V_\bullet,
\]
a holomorphic subbundle of $F\times E$ with fibre $V_j$
over $V_\bullet$.
Their successive quotients
\[
 L_j=\mathcal S_j/\mathcal S_{j-1},\qquad \mathcal S_0=0,
\]
are holomorphic line bundles with fibres $V_j/V_{j-1}$.

We equip $L_j$ with the quotient metric induced by the
standard Hermitian metric on $E$. 
For $\kappa\in K$, write
$\kappa L=(z_1,\ldots,z_\ell)\in K/L\cong(S^1)^\ell$
and define the coordinate characters by
$\chi_j(\kappa L)=z_j$. By the block-diagonal form of $K$,
\[
 \kappa a_j=z_j a_j.
\]
In $V_j^0/V_{j-1}^0$, we have
\[
 [\kappa a_j]=[a_j]\chi_j(\kappa L).
\]
Hence the maps
\[
 \begin{aligned}
 (G/L)\times_{\chi_j}\C&\longrightarrow L_j,\\
 [gL,c]&\longmapsto
 \bigl(gV_\bullet^0,[ga_j]c\bigr)
 \end{aligned}
\]
are well-defined.
For a local smooth section $g:U\to G$ of $G\to F$, the vectors $[g(b)a_j]$ form unit frames of $L_j|_U$. In these frames, the bundle maps are $(b,c)\mapsto(b,c)$, so they are unitary smooth line-bundle isomorphisms.

The resulting unit frames define the map
\[
 \begin{aligned}
 \Phi:G/L&\longrightarrow
 S(L_1)\times_F\cdots\times_F S(L_\ell),\\
 gL&\longmapsto
 \bigl(gV_\bullet^0;[ga_1],\ldots,[ga_\ell]\bigr).
 \end{aligned}
\]
Under the given local trivializations, $\Phi$ has the form
\[
 (b,\kappa L)\longmapsto
 \bigl(b,\chi_1(\kappa L),\ldots,\chi_\ell(\kappa L)\bigr).
\]
Since $(\chi_1,\ldots,\chi_\ell):K/L\to(S^1)^\ell$
is a Lie group isomorphism, $\Phi$ is a diffeomorphism of principal
torus bundles. Thus the holomorphic tautological quotient line bundles realize $G/L$ as the required unit-circle fibre product.


\section{Lens spaces and moment-angle manifolds}\label{ex:nonhomogeneous}

We prove \cref{cor:intro-lens-products,cor:intro-moment-angle} using holomorphic line bundles over products of projective spaces and toric varieties, respectively. In each case, we identify the unit-circle fibre product required by \cref{thm:main}. The examples include products of odd-dimensional spheres and certain connected sums of sphere products, giving explicit smooth factors for the products with $S^4$.

\subsection{Lens spaces and odd-dimensional spheres}

For integers \(a\ge0\) and \(p\ge1\), we define the lens space by
\[
 L^{2a+1}(p)=S^{2a+1}/\Z_p.
\]
Here the group \(\Z_p\) is generated by a primitive \(p\)-th root of unity. It acts by scalar multiplication
on the unit sphere \(S^{2a+1} \subset\C^{a+1}\).
In particular, \(L^{2a+1}(1)=S^{2a+1}\).

\begin{restatedcorollary}{cor:intro-lens-products}
Let $\ell\ge2$ be even, and let $a_j\ge0, p_j\ge1$
be integers for $1\leq j\leq\ell$. Then
\[
 S^4\times\prod_{j=1}^{\ell}L^{2a_j+1}(p_j)
\]
admits a complex structure.
\end{restatedcorollary}

\begin{proof}

We first realize each lens space as a unit circle bundle. We then assemble these bundles over a product of projective spaces. For any $\ell\ge1$, let
\[
 B=\prod_{j=1}^{\ell}\mathbb P^{a_j},
\]
with projections $\operatorname{pr}_j:B\to\mathbb P^{a_j}$.
Pulling back the holomorphic line bundles $\mathcal O_{\mathbb P^{a_j}}(-p_j)$ along
$\operatorname{pr}_j$ defines holomorphic line bundles
\[
 L_j=\operatorname{pr}_j^*
 \mathcal O_{\mathbb P^{a_j}}(-p_j)
\]
over $B$.
Using the standard Hermitian metric on $\C^{a+1}$ and the induced tensor product metric, we have
\[
 S(\mathcal O_{\mathbb P^a}(-1))
 =\{(\lambda,v):\lambda\in\mathbb P^a,\ v\in\lambda,\
 |v|=1\}
 \cong S^{2a+1},
\]
and
\[
 S(\mathcal O_{\mathbb P^a}(-p))
 =\{(\lambda,\xi):\lambda\in\mathbb P^a,\
 \xi\in\lambda^{\otimes p},\ |\xi|=1\}.
\]
The tensor power map
\[
 S^{2a+1}\longrightarrow S(\mathcal O_{\mathbb P^a}(-p)),
 \qquad v\longmapsto(\C v,v^{\otimes p}),
\]
is well defined since $|v^{\otimes p}|=|v|^p$.
We choose a local smooth unit frame
$e(b)$ of $\cO_{\mathbb P^a}(-1)$ and write \(v=z\,e(b)\). Then we have \(v^{\otimes p}=z^p e(b)^{\otimes p}\) with \(e(b)^{\otimes p}\) the unit frame of \(\cO_{\mathbb P^a}(-p) \).
Thus the map is locally 
\[(b,z)\mapsto(b,z^p).\] 
It is therefore a smooth covering whose
fibres are exactly the scalar $\Z_p$-orbits, where $\Z_p$
is the group of $p$th roots of unity. Hence it descends to
a diffeomorphism
\[
 L^{2a+1}(p)=S^{2a+1}/\Z_p
 \cong_{\diff}S(\mathcal O_{\mathbb P^a}(-p)).
\]

Since the $L_j$ are pulled back from the separate factors of $B$, these circle-bundle identifications give
\[
 S(L_1)\times_B\cdots\times_B S(L_\ell)
 \cong_{\diff}
 \prod_{j=1}^{\ell}S(\mathcal O_{\mathbb P^{a_j}}(-p_j))
 \cong_{\diff}
 \prod_{j=1}^{\ell}L^{2a_j+1}(p_j).
\]
There are $\ell$ circle factors. For even $\ell$, \cref{thm:main} applies to this fibre product and proves the corollary.
\end{proof}

Taking all $p_j=1$ recovers sphere factors. In particular, $p_1=p_2=1$, $a_1=1$ and $a_2=3$ give a complex structure on $S^4\times S^3\times S^7$.
Qian--Tang--Yan \cite[Theorems~1.2 and~1.3]{QianTangYan2023}
identify it with the indefinite OT--FKM isoparametric hypersurfaces $M\subset S^{15}$ with multiplicities $(4,3)$.


\subsection{Moment-angle manifolds}

Cox \cite{Cox1995} constructs a holomorphic principal bundle over the toric variety associated with a smooth complete rational fan. Panov--Ustinovsky \cite{PanovUstinovsky2012} identify its compact reduction with a moment-angle manifold. We realize this reduction as a unit-circle fibre product of holomorphic character line bundles and apply \cref{thm:main}. The polygon descriptions of McGavran \cite{McGavran1979} and selected examples of Bosio--Meersseman \cite{BosioMeersseman2006} then give explicit smooth factors for the products with $S^4$.

\begin{restatedcorollary}{cor:intro-moment-angle}
Let $\Sigma$ be a smooth complete rational fan in $\R^d$
with $m$ rays, where $m-d$ is positive and even.
Then $S^4\times\mathcal Z_\Sigma$ admits a complex structure,
with the smooth structure specified in \cref{tor:cox-construction}.
In particular, this applies to $S^4\times\mathcal Z_{P_m}$
for a Delzant $m$-gon $P_m$ with even $m\ge4$, and to
$S^4\times\mathcal Z_P$ when $P$ is obtained from a Delzant $d$-simplex by an odd number $q\ge1$ of vertex cuts, $d\ge2$;
see \eqref{tor:product_polygon} and \eqref{tor:product_truncated-simplex}.
\end{restatedcorollary}

\subsubsection{Character line bundles of the Cox bundle}\label{tor:cox-construction}

Let $\Sigma$ be a smooth complete rational fan in $\R^d$,
with primitive ray generators $a_1,\ldots,a_m\in\Z^d$.
The map
\[
 A:\Z^m\longrightarrow\Z^d,\qquad e_j\longmapsto a_j,
\]
is onto, since the generators of a maximal cone form an
integral basis. With $r=m-d$, the induced torus homomorphisms
give
\[
 G=\ker\bigl((\C^*)^m\to(\C^*)^d\bigr)\cong(\C^*)^r,
 \qquad K=G\cap\T^m\cong\T^r.
\]
The group $G$ acts by coordinatewise multiplication on
the Cox open set
\[
 U_\Sigma=\bigcup_{\sigma\in\Sigma}
 \{z\in\C^m:z_j\ne0\text{ whenever }a_j\notin\sigma\}.
\]
Cox's quotient theorem identifies
$B_\Sigma=U_\Sigma/G$ with the toric variety of $\Sigma$
\cite[Theorem~2.1]{Cox1995}.
Since the fan is smooth,
\[
 G\longrightarrow U_\Sigma\longrightarrow B_\Sigma
\]
is a holomorphic principal bundle. The base $B_\Sigma$ is compact because the fan is complete
\cite[Section~5]{PanovUstinovsky2012}.

Put
\[
 A_\Sigma=\exp\bigl(\ker(A\otimes\R)\bigr)
 \subset(\R_{>0})^m.
\]
By the polar decomposition $G=K\times A_\Sigma$, quotienting by $A_\Sigma$ gives the compact reduction as a smooth principal $K$-bundle
\[
 K\longrightarrow
 \mathcal Z_\Sigma:=U_\Sigma/A_\Sigma
 \longrightarrow B_\Sigma.
\]
We use $\mathcal Z_\Sigma$ as a smooth model of the
moment-angle manifold associated with $\Sigma$
\cite[Theorem~2.2]{PanovUstinovsky2012}.

To apply the main theorem, we identify this compact reduction with the unit-circle fibre product of holomorphic character line bundles. A basis $\chi_1,\ldots,\chi_r$ of the character lattice of $G$ defines a holomorphic Lie group isomorphism
\[
 \Psi:G\longrightarrow(\C^*)^r,\qquad
 g\longmapsto(\chi_1(g),\ldots,\chi_r(g)).
\]
These characters determine holomorphic line bundles
\[
 L_j=U_\Sigma\times_{\chi_j}\C,\qquad
 (qg,v)\sim(q,\chi_j(g)v).
\]
$\Psi$ takes $K$ to $(S^1)^r$ and $A_\Sigma$ to
$(\R_{>0})^r$. For smooth Hermitian metrics on the $L_j$, the map
\[
 \begin{aligned}
 \Phi:\mathcal Z_\Sigma&\longrightarrow
 S(L_1)\times_{B_\Sigma}\cdots\times_{B_\Sigma}S(L_r),\\
 [q]_{A_\Sigma}&\longmapsto
 \left(\frac{[q,1]}{|[q,1]|_{L_1}},\ldots,
       \frac{[q,1]}{|[q,1]|_{L_r}}\right)
 \end{aligned}
\]
is well defined: replacing $q$ by $qa$, with $a\in A_\Sigma$,
multiplies the $j$-th vector and its norm by the same positive
number $\chi_j(a)$. Using a local section of $U_\Sigma\to B_\Sigma$ and the induced unit frames, $\Phi$ has the local expression
\[
 (b,k)\longmapsto(b,\Psi(k)),\qquad k\in K.
\]
Since $\Psi|_K:K\to(S^1)^r$ is a Lie group isomorphism,
$\Phi$ is a diffeomorphism.

\begin{proof}[Proof of \cref{cor:intro-moment-angle} for $\mathcal Z_\Sigma$]

For positive even $r=m-d$, the identification above realizes $\mathcal Z_\Sigma$ as the unit-circle fibre product of $r$ holomorphic line bundles over the compact connected complex base $B_\Sigma$. Applying \cref{thm:main} proves the assertion.

\end{proof}

\subsubsection{Polygons and an explicit hexagon}

Let $P_m$ be a Delzant polygon \cite{Delzant1988} with $m\geq4$ edges, $\Sigma$ be its normal fan and $\mathcal Z_{P_m}$ be the standard smooth moment-angle manifold of $P_m$.
By \cite[Section~4]{PanovUstinovsky2012}, the corresponding moment-angle manifolds are diffeomorphic:
\[
 \mathcal Z_{P_m}\cong_{\diff}
 \mathcal Z_\Sigma=U_\Sigma/A_\Sigma.
\]

\begin{proof}[Proof of \cref{cor:intro-moment-angle} for $\mathcal Z_{P_m}$]

Since the normal fan $\Sigma$ of a Delzant polygon  is smooth, complete, and rational,
with $m$ rays in $\R^2$ and $d=2$, the preceding proof shows that $S^4 \times \mathcal Z_{\Sigma} $ admits a complex structure when $m$ is even.
It follows that $S^4\times\mathcal Z_{P_m}$
admits a complex structure by the diffeomorphism above.

To write the smooth factors in these products explicitly, we use McGavran's diffeomorphism formula \cite[Theorem~3.4]{McGavran1979}
\begin{equation}\label{tor:polygon-sum}
 \mathcal Z_{P_m}\cong_{\diff}
 \mathop{\#}_{j=3}^{m-1}
 (j-2)\binom{m-2}{j-1}
 (S^j\times S^{m+2-j}).
\end{equation}
Therefore, for every even integer $m\ge4$, the smooth manifold
\begin{equation}\label{tor:product_polygon}
 S^4\times
 \left(
 \mathop{\#}_{j=3}^{m-1}
 (j-2)\binom{m-2}{j-1}
 (S^j\times S^{m+2-j})
 \right)
\end{equation}
admits a complex structure.
\end{proof}

For a Delzant hexagon, we make the Cox presentation and its four character line bundles explicit. This gives the $m=6$ case directly through the main theorem.

\begin{figure}[htbp]
\centering
\begin{minipage}[t]{0.48\textwidth}
\centering
\begin{tikzpicture}[
  scale=1.45,
  >=stealth,
  every node/.style={font=\small}
]
\path[use as bounding box]
  (-1.65,-1.55) rectangle (1.65,1.45);

\draw[thick]
  (-1,0)--(-1,1)--(0,1)--(1,0)
  --(1,-1)--(0,-1)--cycle;
\node at (0,0) {$P$};

\node[left]        at (-1,.5)   {$F_1$};
\node[below left]  at (-.5,-.5) {$F_2$};
\node[below]       at (.5,-1)   {$F_3$};
\node[right]       at (1,-.5)   {$F_4$};
\node[above right] at (.5,.5)   {$F_5$};
\node[above]       at (-.5,1)   {$F_6$};

\draw[->,thick] (-1,.5)--(-.65,.5)
  node[midway,above] {$a_1$};
\draw[->,thick] (-.5,-.5)--(-.25,-.25)
  node[midway,below right] {$a_2$};
\draw[->,thick] (.5,-1)--(.5,-.65)
  node[midway,right] {$a_3$};
\draw[->,thick] (1,-.5)--(.65,-.5)
  node[midway,above] {$a_4$};
\draw[->,thick] (.5,.5)--(.25,.25)
  node[midway,above left] {$a_5$};
\draw[->,thick] (-.5,1)--(-.5,.65)
  node[midway,right] {$a_6$};

\node at (0,-1.45) {(a) Delzant hexagon $P$};
\end{tikzpicture}
\end{minipage}
\hfill
\begin{minipage}[t]{0.48\textwidth}
\centering
\begin{tikzpicture}[
  scale=1.1,
  >=stealth,
  every node/.style={font=\small}
]
\path[use as bounding box]
  (-2.6,-2.05) rectangle (2.6,1.9);

\draw[->,thick] (0,0)--(1.65,0)
  node[right] {$a_1=(1,0)$};
\draw[->,thick] (0,0)--(1.25,1.25)
  node[above right] {$a_2=(1,1)$};
\draw[->,thick] (0,0)--(0,1.65)
  node[left] {$a_3=(0,1)$};
\draw[->,thick] (0,0)--(-1.65,0)
  node[left] {$a_4=(-1,0)$};
\draw[->,thick] (0,0)--(-1.25,-1.25)
  node[below left] {$a_5=(-1,-1)$};
\draw[->,thick] (0,0)--(0,-1.65)
  node[right] {$a_6=(0,-1)$};

\foreach \x/\y in
  {1/0,1/1,0/1,-1/0,-1/-1,0/-1}
  \fill (\x,\y) circle (1.3pt);

\fill (0,0) circle (1.3pt);
\node[below right] at (0,0) {$0$};

\node at (.95,.38)   {$\sigma_{12}$};
\node at (.38,.95)   {$\sigma_{23}$};
\node at (-.75,.75)  {$\sigma_{34}$};
\node at (-.95,-.38) {$\sigma_{45}$};
\node at (-.38,-.95) {$\sigma_{56}$};
\node at (.75,-.75)  {$\sigma_{61}$};

\node at (0,-2.0) {(b) normal fan};
\end{tikzpicture}
\end{minipage}

\end{figure}
The hexagon is
$P=\{(x,y)\in\mathbb R^2:
|x|\le1,\ |y|\le1,\ |x+y|\le1\}$, with $a_i$ inward normals.
The edge $F_i$ and the vertex
$F_i\cap F_{i+1}$ correspond to the ray
$\mathbb R_{\ge0}a_i$ and cone
$\sigma_{i,i+1}=\operatorname{cone}(a_i,a_{i+1})$ in the normal fan, respectively. The ray generators lie in $\Z^2$, so the fan is rational. These cones cover $\R^2$, and each consecutive
pair has determinant one. The fan is therefore complete
and smooth. 

The ray generators determine the matrix
\[
 A=
 \begin{pmatrix}
  1&1&0&-1&-1&0\\
  0&1&1&0&-1&-1
 \end{pmatrix}
\]
and the induced homomorphism 
\[
 A_{\C}:(\C^*)^6\longrightarrow(\C^*)^2,\qquad
 (g_1,\ldots,g_6)\longmapsto
 \left(\frac{g_1g_2}{g_4g_5},
       \frac{g_2g_3}{g_5g_6}\right).
\]
$G=\ker A_{\C}$ is defined by
\[
 g_1=g_2^{-1}g_4g_5,\qquad
 g_3=g_2^{-1}g_5g_6.
\]
Therefore we can write the complex Lie group isomorphism
\[
 \begin{aligned}
 (\C^*)^4&\longrightarrow G,\\
 (h_1,h_2,h_3,h_4)&\longmapsto
 (h_1^{-1}h_2h_3,\ h_1,\ h_1^{-1}h_3h_4,\
  h_2,\ h_3,\ h_4).
 \end{aligned}
\]
For this fan, the Cox open set is
\[
 U_\Sigma=
 \bigcup_{i=1}^{6}
 \{z\in\C^6:z_j\ne0\text{ for }j\notin\{i,i+1\}\},
\]
where indices are read cyclically modulo $6$.
The right action of $G$ on $U_\Sigma$ is, for \(h=(h_1,\ldots,h_4)\in(\C^*)^4\),
\[
 \begin{aligned}
 z\cdot h={}&
 (h_1^{-1}h_2h_3z_1,\ h_1z_2,\
 h_1^{-1}h_3h_4z_3,h_2z_4,\ h_3z_5,\ h_4z_6).
 \end{aligned}
\]

The coordinate characters $\chi_j(h)=h_j$ induce the four holomorphic line bundles to which we apply the main theorem
\[
 L_j=U_\Sigma\times_{\chi_j}\C,\qquad 1\le j\le4.
\]
By the proof of \cref{cor:intro-moment-angle},
\[
 \mathcal Z_{P_6}\cong_{\diff}
 S(L_1)\times_{B_\Sigma}\cdots\times_{B_\Sigma}S(L_4).
\]
Applying \cref{thm:main} to these four lines and using the diffeomorphism \eqref{tor:polygon-sum}, we obtain a complex structure on
\[
 S^4\times
 \left(
 \#^9(S^3\times S^5)
 \#^8(S^4\times S^4)
 \right).
\]

\subsubsection{Delzant realizations and Bosio--Meersseman examples}
\label{tor:bm-examples}

We now apply the fan construction to two families described by Bosio--Meersseman.
The diffeomorphism type of $\mathcal Z_P$ depends only on the combinatorial type of $P$
\cite[Theorem~4.1 and Corollary~4.7]{BosioMeersseman2006}.
It therefore suffices to choose Delzant realizations of the polytopes involved.
Such realizations exist for the examples below, since products and suitably chosen face truncations preserve the Delzant property; see \cite[Section~5.2]{BuchstaberPanov2015}.
For each realization, the normal fan $\Sigma$ is smooth, complete, and rational, and
$\mathcal Z_P\cong_{\diff}\mathcal Z_\Sigma$ \cite[Section~4]{PanovUstinovsky2012}.
The relevant rank is $r=m-d$, where $m$ is the number of facets and $d=\dim P$.

\begin{enumerate}
\item Products of simplices give products of odd-dimensional spheres:
\[
 P=\Delta^{d_1}\times\cdots\times\Delta^{d_a},\qquad
 \mathcal Z_P\cong_{\diff}\prod_{j=1}^{a}S^{2d_j+1},
 \qquad d_j\geq1.
\]
These include the sphere products in
\cite[Example~0.4 and Corollary~4.6]{BosioMeersseman2006}.
Here $r=a$. For even $a$, the resulting products with $S^4$ recover the preceding lens-space examples with all $p_j=1$.
For odd $a$, \cref{cor:odd-rank} gives complex structures after adding an $S^1$ factor.

\item Let $P$ be obtained from a $d$-simplex by $q\geq1$ successive vertex cuts, where $d\geq2$.
Then $P$ has $m=d+q+1$ facets.
Bosio--Meersseman \cite[Theorem~6.3]{BosioMeersseman2006} extend McGavran's polygon formula to the decomposition
\begin{equation}\label{tor:truncated-simplex}
 \mathcal Z_P\cong_{\diff}
 \mathop{\#}_{j=1}^{q}
 j\binom{q+1}{j+1}
       (S^{j+2}\times S^{2d+q-j-1}).
\end{equation}
The case $d=2$ recovers the polygon family above.
Since $r=q+1$, an odd number of vertex cuts gives the even rank required in \cref{cor:intro-moment-angle}.
\end{enumerate}

\begin{proof}[Proof of \cref{cor:intro-moment-angle} for $\mathcal Z_P$]
For the truncated simplices in (2), $r=q+1$ is positive and even when $q\geq1$ is odd.
The fan case therefore gives a complex structure on $S^4\times\mathcal Z_P$.
Using \eqref{tor:truncated-simplex}, we obtain complex structures on
\begin{equation}\label{tor:product_truncated-simplex}
    S^4 \times \left(\mathop{\#}_{j=1}^{q}
 j\binom{q+1}{j+1}
       (S^{j+2}\times S^{2d+q-j-1})\right).
\end{equation}
This proves the remaining assertion of \cref{cor:intro-moment-angle}.
\end{proof}

We next compare the homeomorphism types of the manifolds
in \cref{cor:intro-moment-angle} with those of the complex
manifolds constructed by Bosio and Meersseman \cite{BosioMeersseman2006}.
Let $\Sigma$ be a smooth complete rational fan with
$r=m-d>0$, and put $N=S^4\times\mathcal Z_\Sigma$.
Our preceding construction equips $M=N$ with a complex
structure when $r$ is even, and $M=N\times S^1$ when
$r$ is odd.

Bosio--Meersseman construct complex structures on
even-dimensional links and, for an odd-dimensional link $W$,
on $W\times S^1$ and $W/S^1$
\cite[Theorem~12.2]{BosioMeersseman2006}.
Every link is diffeomorphic to $\T^\ell\times\mathcal Z_Q$,
where $Q$ is a simple convex polytope,
$\mathcal Z_Q$ is its standard moment-angle manifold,
and $\ell\ge0$.

Both $\mathcal Z_\Sigma$ and $\mathcal Z_Q$ are $2$-connected
\cite[Proposition~4.3.5(a)]{BuchstaberPanov2015}.
The Hurewicz and universal coefficient theorems, together
with the K\"unneth formula, imply $H^2(M;\Z)=0$.
Moreover, the universal covers of $M$ and
$\T^\ell\times\mathcal Z_Q$ are homotopy equivalent to
$N$ and $\mathcal Z_Q$, respectively.

After tensoring with $\mathbb Q$, the decomposition and
product rule in \cite[Theorem~4.5.8]{BuchstaberPanov2015}
imply that $H^4(\mathcal Z_Q;\mathbb Q)$ is spanned by
classes annihilating nonzero degree-three classes.
Indeed, each nonzero degree-four summand corresponds to
a disconnected full subcomplex on three vertices.
Two of these vertices form a missing edge, which determines
a nonzero degree-three class. Its product with every class
in the summand vanishes by the product rule, since the
two vertex sets intersect.
However, this fails for $N$.
Write a degree-four class as $cu+v$ in $H^4(N;\mathbb Q)$, where $u$ is pulled back from a generator of $H^4(S^4;\mathbb Q)$ and $v$
from $\mathcal Z_\Sigma$.
If $(cu+v) \smile a=0$ for some nonzero $a\in H^3(N;\mathbb Q)$,
the K\"unneth formula forces $c=0$.
Thus these classes cannot span $H^4(N;\mathbb Q)$.
$N$ and $\mathcal Z_Q$ have nonisomorphic graded
rational cohomology rings.

A homeomorphism from $M$ to a link, or to its product
with $S^1$, would contradict this obstruction upon
passing to universal covers.
If $M$ were homeomorphic to $W/S^1$, the principal
circle bundle $W\to W/S^1$ would pull back to a trivial
bundle over $M$, since $H^2(M;\Z)=0$. So $W\cong_{\mathrm{homeo}}M\times S^1$, and their universal
covers would give the same contradiction.
Therefore, $M$ is not homeomorphic to any manifold in
the Bosio--Meersseman construction.

\vspace{2mm}
\needspace{4\baselineskip}
\noindent
\textbf{AI disclosure.}
The authors used AI tools to develop the complex structure on $S^4 \times V_{2s}(\HH^n)$ during extensive discussions on possible generalizations of their fundamental work on $S^4 \times \T^2$ \cite{GaoTangYan2026}, and check the later examples of the general construction formulated by the authors. The authors checked and refined all arguments and take full responsibility for the mathematical content.

\needspace{5\baselineskip}


\end{document}